\documentclass[a4paper, 12pt, reqno]{amsart}
\usepackage[top=34truemm,bottom=30truemm,left=30truemm,right=30truemm]{geometry}

\usepackage{amscd, amsmath, amssymb, amsthm, url, hyperref}

\numberwithin{equation}{section}

\theoremstyle{plain}
\newtheorem{theorem}{Theorem}[section]
\newtheorem{proposition}[theorem]{Proposition}
\newtheorem{corollary}[theorem]{Corollary}
\newtheorem{lemma}[theorem]{Lemma}

\theoremstyle{definition}
\newtheorem{definition}[theorem]{Definition}

\newtheorem{remark}[theorem]{Remark}
\newtheorem*{rem*}{Remark}

\theoremstyle{remark}

\begin{document}

\title[No exceptional orbits under polar actions on Hilbert spaces]{Non-existence of exceptional orbits \\ under polar actions on Hilbert spaces}
\author[M. Morimoto]{Masahiro Morimoto}



\thanks{Supported by the Grant-in-Aid for Research Activity Start-up (No.\ 20K22309), by Grant-in-Aid for JSPS Research fellow (No.\ 23KJ1793) and by MEXT Promotion of Distinctive Joint Research Center Program JPMXP0723833165}





\maketitle

\ \vspace{-7mm}
\begin{center}
\footnotesize
Department of Mathematical Sciences, Tokyo Metropolitan University
\\
1-1 Minami-Osawa, Hachioji-shi, Tokyo 192-0397, Japan
\\
email: morimoto.mshr@gmail.com 
\vspace{-1mm}
\end{center}

 \smallskip

\begin{abstract}
We prove that any polar action on a separable Hilbert space by a connected Hilbert Lie group does not have exceptional orbits. This generalizes a result of Berndt, Console and Olmos in the finite dimensional Euclidean case. As an application, we give a simple geometric proof of the fact that any hyperpolar action on a simply connected compact Riemannian symmetric space by a connected Lie group does not have exceptional orbits.
\end{abstract}

 \medskip \smallskip

\begin{center}
\footnotesize
\emph{keywords}: polar action, Hilbert space, exceptional orbit, isoparametric submanifold
\bigskip
\\
\emph{2020 Mathematics Subject Classification}: 53C40.
\end{center}

 \medskip

\section{Introduction}

A proper isometric action of a Lie group on a complete Riemannian manifold $M$ is called \emph{polar} if there exists a connected complete immersed submanifold $\Sigma$ of $M$ which meets every orbit and is orthogonal to the orbits at every point of intersection. Such a $\Sigma$ is called a \emph{section}, which is automatically totally geodesic in $M$. If $\Sigma$ is flat with respect to the induced metric, then the action is called \emph{hyperpolar}. For details of polar actions and hyperpolar actions, see \cite{Gor22} and its references.

In 2003, Berndt, Console and Olmos showed the following fact as one of the important properties of polar actions on Euclidean spaces \cite[Remark 3.2.10 and Corollary 5.4.3]{BCOI} (see also the second edition \cite[Remark 2.3.12 and Corollary 4.4.3]{BCO}):


\begin{theorem}[Berndt-Console-Olmos \cite{BCOI, BCO}]\label{thmBCO}
Let $G$ be a connected Lie group acting properly and isometrically on the Euclidean space $\mathbb{E}^n$. If the action of $G$ on $\mathbb{E}^n$ is polar, then it has no exceptional orbits.
\end{theorem}


Here an \emph{exceptional} orbit is by definition a maximal dimensional orbit which is not principal. Note that if $G$ is not connected, there exists a counterexample to the theorem (see Remark \ref{rembco}). Their proof of Theorem \ref{thmBCO} is based on the theory of isoparametric submanifolds in Euclidean spaces investigated by Terng \cite{Ter85}. Here a submanifold of $\mathbb{E}^n$ is called \emph{isoparametric} if the normal bundle is flat and the principal curvatures in the direction of any parallel normal vector field are constant. They showed that any maximal dimensional orbit of a polar action on $\mathbb{E}^n$ is an isoparametric submanifold of $\mathbb{E}^n$  and proved Theorem \ref{thmBCO} by using Terng's result that any connected isoparametric submanifold of $\mathbb{E}^n$ has a globally flat normal bundle. 

The theory of polar actions on Euclidean spaces together with the theory of isoparametric submanifolds in Euclidean spaces was generalized to the case of Hilbert spaces by Palais and Terng \cite{PT88, Ter89, Ter95}. They introduced a suitable class of Lie group actions on Hilbert manifolds, namely \emph{proper Fredholm} (PF) actions, and a suitable class of submanifolds in Hilbert spaces, namely \emph{proper Fredholm} (PF) submanifolds. 
It follows that every orbit of an isometric PF action on a separable Hilbert space $V$ is a PF submanifold of $V$. They investigated polar PF actions on Hilbert spaces and isoparametric PF submanifolds in Hilbert spaces and gave many interesting examples of them. Then it is a natural question whether Theorem \ref{thmBCO} can be generalized to the case of polar PF actions on Hilbert spaces. The main purpose of this paper is to give a positive answer to this question. 

To explain our result, we review basic facts about PF actions \cite[Section 5]{PT88}. A smooth action of a Hilbert Lie group $\mathcal{L}$ on a Hilbert manifold $M$ is called proper Fredholm (PF) if the map $\mathcal{L} \times M \rightarrow M \times M$, $(g, p) \mapsto (g \cdot p, p)$ is proper and if for each $p \in M$ the map $\mathcal{L} \rightarrow M$, $g \mapsto g \cdot p$ is Fredholm. Then any orbit $\mathcal{L} \cdot p$ has finite codimension and the isotropy subgroup $\mathcal{L}_p$ at $p$ is compact. An orbit $\mathcal{L} \cdot p$ is called \emph{principal} if there exists a neighborhood $U$ of $p$ such that for all $q \in U$ there exists $g \in \mathcal{L}$ such that $\mathcal{L}_p \subset g \mathcal{L}_q g^{-1}$. Although this definition is local, we can assume $U = M$ when $M$ is connected (Lemma \ref{lem0}). A point $p \in M$ is called \emph{principal} if the orbit $\mathcal{L} \cdot p$ is principal. If the $\mathcal{L}$-action on $M$ is PF, then it follows that there exists a slice at each point of $M$, there exists a principal orbit type of the $\mathcal{L}$-action, and the set of all principal points is open and dense in $M$. In particular, if $M$ is connected, then the principal orbit type is unique.

To define exceptional orbits of PF actions, we remark that an orbit of ``maximal dimension"  is not meaningful for PF actions in general. Thus we consider an orbit of ``minimal codimension" instead. We assume that $M$ is connected and prove that principal orbits of PF actions have minimal codimension (Proposition \ref{cor0}). Based on this fact, we define an exceptional orbit as an orbit of minimal codimension which is not principal (Definition \ref{def0}). Then we prove the following theorem which generalizes Theorem \ref{thmBCO} to the case of polar PF actions on Hilbert spaces. Here, an isometric PF action on a separable Hilbert space $V$ is called \emph{polar} if there exists a closed affine subspace of $V$ which meets every orbit orthogonally.


\begin{theorem}\label{thm2}
Let $\mathcal{L}$ be a connected Hilbert Lie group acting isometrically and PF on a separable Hilbert space $V$.  If the action of $\mathcal{L}$ on $V$ is polar, then it has no exceptional orbits.
\end{theorem}


Note that if $V$ is finite dimensional, then the Fredholm condition is automatically satisfied and Theorem \ref{thm2} is the same as Theorem \ref{thmBCO}. Note also that if $\mathcal{L}$ is not connected, there exist counterexamples (Remarks \ref{rembco} and \ref{rem5}). We prove Theorem \ref{thm2} based on the theory of isoparametric PF submanifolds in Hilbert spaces investigated by Palais and Terng \cite{PT88, Ter89} and also by Heintze, Liu and Olmos \cite{HLO06}. The strategy of proof of Theorem \ref{thm2} is similar to that of Theorem \ref{thmBCO}. We make clear some inaccuracies in the original proof and give a complete proof in the case of Hilbert space  (see Remarks \ref{rem8} and \ref{rembco2}). Our proof is valid also in the finite dimensional Euclidean case. 

As an application of Theorem \ref{thm2}, we give a simple geometric proof of the fact that hyperpolar actions on simply connected compact Riemannian symmetric spaces by connected Lie groups does not have exceptional orbits. Recall that Terng \cite{Ter95} associated to each hyperpolar action of a compact Lie group $H$ on a compact symmetric space $G/K$, a polar PF action of an infinite dimensional Hilbert Lie group $P(G, H \times K)$ on the Hilbert space $V_\mathfrak{g} = L^2([0,1], \mathfrak{g})$ consisting of all $L^2$-paths with values in the Lie algebra $\mathfrak{g}$ of $G$ (see also \cite{PT88, Ter89, PiTh90}). It follows that an $H$-orbit is exceptional if and only if the corresponding $P(G, H \times K)$-orbit is exceptional. We prove that $P(G, H \times K)$ is connected if $G$ is simply connected and if $K$ and $H$ are connected (Lemma \ref{lemconnected}). Thus, by applying Theorem \ref{thm2} to the $P(G, H \times K)$-action, we obtain:


\begin{corollary}\label{thm1}
Let $G$ be a simply connected compact Lie group and $K$ a connected symmetric subgroup of $G$ with involution $\theta$. Assume that $G$ is equipped with a bi-invariant Riemannian metric induced from an inner product of $\mathfrak{g}$ invariant by $\operatorname{Ad}(G)$ and $\theta$ (and thus $G/K$ becomes a simply connected compact Riemannian symmetric space). Let $H$ be a closed connected subgroup of $G$. If the  action of $H$ on $G/K$ is hyperpolar, then it has no exceptional orbits. 
\end{corollary}


We now mention the related results by other researchers. 

In 2003, Gorodski and Thorbergsson \cite[Proposition 3.3]{GT03} proved that a taut representation of a connected compact Lie group on Euclidean spaces does not have exceptional orbits. This gives another proof of Theorem \ref{thmBCO} because isoparametric submanifolds in Euclidean spaces are taut \cite{HPT88}. It is possible to generalize their method to the infinite dimensional case because the square of the distance function from a point to a PF submanifold satisfies the Palais-Smale condition and thus the Morse theory also holds on infinite dimensional PF submanifolds \cite{PT88, Ter89} (see also \cite[Section 5]{M23}). To prove Theorem \ref{thm2} by their method, one also needs a fact that isoparametric PF submanifolds in Hilbert spaces are taut \cite{PT88, Ter89} (see also \cite[Appendix B]{HLO06} for a remark in the case of non-injective immersion). In our proof of Theorem \ref{thm2}, we do not use such a fact or Morse theory at all. The only nontrivial fact we use is the global flatness of normal bundles of isoparametric PF submanifolds in Hilbert spaces \cite{HLO06}. 

Let $M$ be a finite dimensional \emph{simply connected} complete Riemannian manifold. In 2006, Alexandrino and T\"oben \cite[Theorem 1.5]{AT06} studied singular Riemannian foliations with sections and proved that if $\mathcal{F}$ is such a foliation on $M$ and if the leaves of $\mathcal{F}$ are compact, then each maximal dimensional leaf of $\mathcal{F}$ has a globally flat normal bundle. In 2011, Alexandrino \cite[Corollary 1.3]{Ale11} gave another proof of this result without assuming that the leaves are compact. Using his result, one can prove the non-existence of exceptional orbits under polar actions on $M$ by connected Lie groups, by the arguments similar to those in the proof of Theorem \ref{thm2}. In 2012, in the study of construction of polar actions, Grove and Ziller \cite[Theorem 3.4]{GZ12} showed that polar actions on $M$ by connected Lie groups does not have exceptional orbits. It is not clear whether their methods can be generalized to the infinite dimensional case. Although Corollary \ref{thm1} is contained in their result, it gives a simple geometric proof in the case of hyperpolar actions on compact symmetric spaces.

In this paper, we define exceptional orbits of PF actions and show their properties in Section \ref{sec0}. Then we prove Theorem \ref{thm2} in Section \ref{secII} and Corollary \ref{thm1} in Section \ref{secI}.

\section{Exceptional orbits of  PF actions}\label{sec0}

Throughout this section, $\mathcal{L}$ denotes a Hilbert Lie group acting on a \emph{connected} Hilbert manifold $M$ and the action is supposed to be PF. 
We assume the basic facts about PF actions shown in \cite[Section 5]{PT88}.

First we remark the definition of principal orbits. According to \cite[Definition 5.1.10]{PT88} the orbit $\mathcal{L} \cdot p$ through $p \in M$ is called \emph{principal} if it satisfies the following condition: 
\begin{itemize}
\item[(a)] there exits a neighborhood $U$ of $p$ such that for all $q \in U$ there exists $g \in \mathcal{L}$ such that $\mathcal{L}_p \subset g \mathcal{L}_q g^{-1}$.
\end{itemize}
On the other hand, in the standard definition, the following condition is required: 
\begin{itemize}
\item[(b)] for all $q \in M$ there exists $g \in \mathcal{L}$ such that $\mathcal{L}_p \subset g \mathcal{L}_q g^{-1}$.
\end{itemize}
We prove:

\begin{lemma}\label{lem0}
\textup{(a)} and \textup{(b)} are equivalent.
\end{lemma}

\begin{proof}
It suffices to show that (a) implies (b). Let $q \in M$ and take a slice $S$ at $q$ \cite[Theorem 5.2.6]{PT88}. Take $s \in S$ such that the orbit $\mathcal{L} \cdot s$ is principal. From the property of slice \cite[Corollary 5.1.13]{PT88} we have $\mathcal{L}_s \subset \mathcal{L}_q$. From the uniqueness of principal orbit type on connected manifolds \cite[Corollary 5.4.18]{PT88} it follows that there exits $g \in \mathcal{L}$ such that $\mathcal{L}_p = g \mathcal{L}_s g^{-1}$. Thus $\mathcal{L}_p \subset g \mathcal{L}_qg^{-1}$. This proves the lemma.
\end{proof}

To define exceptional orbits, we remark that an orbit of ``maximal dimension" is not meaningful for PF actions. Thus we consider an orbit of ``minimal codimension'' instead. We prove:


\begin{lemma}\label{lem1}
The orbit $N = \mathcal{L} \cdot p$ through $p \in M$ has minimal codimension if and only if the isotropy subgroup $\mathcal{L}_p$ of $\mathcal{L}$ at $p$ has minimal dimension. 
\end{lemma}

\begin{proof}
Take a slice $S_p$ at $p$. Then $S_p$ is $\mathcal{L}_p$-invariant and the orbit space $S_p/\mathcal{L}_p = \mathcal{L} S_p/\mathcal{L}$ is an open subset of $M/\mathcal{L}$ \cite[Corollary 5.1.13]{PT88}. Since $M/\mathcal{L}$ is connected, it follows that $\dim (S_p/\mathcal{L}_p)$ does not depend on $p$. Here the dimension of the orbit space is defined to be the dimension of the manifold consisting of all principal points of the action. In particular $\dim (S_p/\mathcal{L}_p) = \dim S_p - \dim \mathcal{L}_p + \dim H$, where $H$ denotes the principal isotropy class of the $\mathcal{L}_p$-action on $S_p$ or equivalently of the $\mathcal{L}$-action on $M$. Thus $\operatorname{codim} N (= \dim S_p)$ is minimal if and only if $\dim \mathcal{L}_p$ is minimal. 
\end{proof}


The following fact is fundamental:


\begin{proposition}\label{cor0}
Principal orbits of the $\mathcal{L}$-action on $M$ have minimal codimension.
\end{proposition}
\begin{proof}
If $\mathcal{L} \cdot p$ is principal then it follows from Lemma \ref{lem0} that for all $q \in M$ there exists $g \in \mathcal{L}$ such that $\mathcal{L}_p \subset g \mathcal{L}_q g^{-1}$. This together with Lemma \ref{lem1} shows that principal orbits have minimal codimension.
\end{proof}

Based on this fact, we make the following definition:


\begin{definition}\label{def0}
An orbit $\mathcal{L} \cdot p$ is called \emph{exceptional} if it has minimal codimension and is not principal.
\end{definition}


Assume that $M$ is a connected Riemannian Hilbert manifold and the action of $\mathcal{L}$ on $M$ is isometric and PF. The representation of $\mathcal{L}_p$ on the normal space $T^\perp_p N$ by differential is called the \emph{slice representation} of the orbit $N = \mathcal{L} \cdot p$. It follows that $N$ is a principal orbit if and only if the slice representation is trivial \cite[Proposition 5.4.7]{PT88}. The following lemma concerns the general case that the orbits have minimal codimension and will be useful later.

\begin{lemma} \label{lem2}
The orbit $N = \mathcal{L} \cdot p$ through $p \in M$ has minimal codimension if and only if  the slice representation of $N$ restricted to the identity component $\mathcal{L}_p^0$ is trivial.
\end{lemma}

\begin{proof}
Take a normal slice $S_p$ at $p$ \cite[Proposition 5.4.3]{PT88}. From the property of slice \cite[Corollary 5.1.13]{PT88}, $\mathcal{L}_q \subset \mathcal{L}_p$ and thus $\mathcal{L}_q^0 \subset \mathcal{L}_p^0$  for all $q \in S_p$. Hence by Lemma \ref{lem1}, $N$ has minimal codimension if and only if $\mathcal{L}_q^0 = \mathcal{L}_p^0$  for all $q \in S_p$. This is equivalent to the condition that $\mathcal{L}_p^0$ fixes every point of $S_p$ and thus to the condition that the $\mathcal{L}_p^0$-representation on $T^\perp_p N$ is trivial. This proves the lemma.
\end{proof}

\section{Proof of Theorem \ref{thm2}} \label{secII}

Throughout this section, $\mathcal{L}$ denotes a connected Hilbert Lie group acting on a separable Hilbert space $V$ and the action is assumed to be isometric and PF.

If $N = \mathcal{L} \cdot p$ is a principal orbit, then every normal vector can be extended to an equivariant normal vector field on $N$ \cite[Corollary 5.4.9]{PT88}. Moreover, if the action is polar, then equivariant normal vector fields on principal orbits are parallel with respect to the normal connection \cite[Theorem 5.6.7]{PT88}. The following lemma concerns the general case that orbits have minimal codimension. This is essential in our proof and will be used twice (see the proofs of Lemma \ref{lem11} and Theorem \ref{thm2} below).


\begin{lemma}\label{lem3}
Let $N = \mathcal{L} \cdot p$ be an orbit of minimal codimension. Then for each $q \in N$ there exists an open neighborhood $U$ of $e \in \mathcal{L}$ such that for each $\xi \in T^\perp_q N$ a normal vector field $\tilde{\xi}$ on the open neighborhood $W := U \cdot q$ in $N$ is well-defined by $\tilde{\xi}(g \cdot q) = dg (\xi)$ where $g \in U$. Moreover, if the $\mathcal{L}$-action is polar, then $\tilde{\xi}$ is parallel with respect to the normal connection.
\end{lemma}


\begin{remark}\label{rem8}
In the proof of \cite[Proposition 2.3.5]{BCO}, it is stated that $\xi$ can be extended to a $G$-equivariant normal vector field $\tilde{\xi}$ on an open neighborhood of $p$ in $G \cdot p$. However $\tilde{\xi}$ is in general not equivariant under the whole $G$ and thus should be stated as in Lemma \ref{lem3}. 
\end{remark}

\begin{proof}[Proof of Lemma \textup{\ref{lem3}}]
Take an open neighborhood $B$ of $e \in \mathcal{L}$ so that $B \cap \mathcal{L}_q \subset \mathcal{L}_q^0$. Take also an open neighborhood  $U$ of $e \in \mathcal{L}$ so that $\{ g^{-1}h \mid g, h \in U \} \subset B$. By Lemma \ref{lem2} it follows that $\tilde{\xi}$ is a well-defined normal vector filed on $W$. Let $r \in W$ and $v \in T^\perp_r N$. Then there exists $X \in \operatorname{Lie}(\mathcal{L})$ such that $v = X^*_r$ where $\operatorname{Lie}(\mathcal{L})$ denotes the Lie algebra of $\mathcal{L}$ and $X^*$ the associated vector field on $V$ defined by
$
X^*_r := 
\left.
\frac{d}{dt}
\right|_{t = 0}
(\exp tX) \cdot r.
$
Clearly $X^*$ is tangent to every $\mathcal{L}$-orbit and thus normal to the section $\Sigma_r$ through $r$. 
Let $\hat{\xi}$ be an extension of $\tilde{\xi}$ to a smooth vector field on $V$. 
Denote by $X^*_t$ the one parameter group of $X^*$. Then
\begin{equation*}
[X^* , \hat{\xi}]_r 
=
\left.
\frac{d}{dt}
\right|_{t =0}
d X^*_{-t} (\hat{\xi}(X^*_t(r)))
=
\left.
\frac{d}{dt}
\right|_{t =0}
d \exp (-tX) (\tilde{\xi}((\exp tX) \cdot r))
=0
\end{equation*}
because for each $g \in U$, $(\exp tX) g \in U$ for sufficiently small $t >0$. Thus we have 
\begin{equation*}
\nabla^\perp_{v} \tilde{\xi}
=
(D_{X^*} \hat{\xi})^\perp_r
=
(D_{\hat{\xi}} X^*)_r^\perp
=
- A^{\Sigma_r}_{X^*(r)} (\hat{\xi}(r))
=
0,
\end{equation*}
where $D$ denotes the flat connection on $V$ and $A^{\Sigma_r}$ the shape operator of $\Sigma_r$. Thus $\tilde{\xi}$ is parallel and the lemma is proved.
\end{proof}

We know that every orbit of an isometric PF action on $V$ is a \emph{proper Fredholm} (PF) submanifold of $V$ \cite[Theorem 7.1.6]{PT88}. Moreover, if the action is polar, then any principal orbit is \emph{isoparametric} \cite{Ter85, PT88, Ter89}, that is, the normal bundle is flat and the principal curvatures of $N$ in the direction of any parallel normal vector field are constant \cite[Theorem 5.7.1]{PT88}. We prove:


\begin{lemma}\label{lem11}
Suppose that the $\mathcal{L}$-action on $V$ is polar. Then an orbit $N = \mathcal{L} \cdot p$ of minimal codimension is an isoparametric PF submanifold of $V$.
\end{lemma}

\begin{proof}
From Lemma \ref{lem3}, every normal vector $\xi$ can be locally extended to a parallel normal vector field $\tilde{\xi}$. This shows that $T^\perp N$ is flat. Moreover, by the local equivariant property of $\tilde{\xi}$, the eigenvalues of the shape operator $A^N_{\tilde{\xi}}$ is constant on the neighborhood $W$. Thus $N$ is isoparametric.
\end{proof}

The following fact was shown by Terng \cite[Proposition 3.6]{Ter85} in the finite dimensional Euclidean case and later generalized by Heintze, Liu and Olmos \cite[Theorem B]{HLO06} to the case of Hilbert space. 


\begin{proposition}[Terng \cite{Ter85}, Heintze-Liu-Olmos \cite{HLO06}]\label{HLO}
Any connected isoparametric PF submanifold of a separable Hilbert space has a globally flat normal bundle.
\end{proposition}


We are now in a position to prove Theorem \ref{thm2}.

\begin{proof}[\textbf{\textup{Proof of Theorem \ref{thm2}}}]
Let $N = \mathcal{L} \cdot p$ be an orbit of minimal codimension. Since $\mathcal{L}$ is connected, $N$ is connected. Then it follows from Lemma \ref{lem11} and Proposition \ref{HLO} that $T^\perp N$ is globally flat. Thus there exist parallel normal frame fields $\{\xi_i\}_{i = 1}^n$ globally defined on $N$. We prove that
\begin{equation*}
\mathcal{L}(\xi_i) : = \{g \in \mathcal{L} \mid \xi_i(g \cdot p) = dg (\xi_i(p))\}
\end{equation*}
is open and closed in $\mathcal{L}$.

Let $g \in \mathcal{L}(\xi_i)$. Set $q = g \cdot p$. By Lemma \ref{lem3} and the uniqueness of the parallel translation there exists an open neighborhood $U$ of $e \in \mathcal{L}$ such that $\xi_i(h \cdot q) = dh (\xi_i (q))$ for all $h \in U$. This together with $g \in \mathcal{L}(\xi_i)$ implies that $\xi_i(h \cdot g \cdot p) = dh \circ dg (\xi_i(p))$. This shows that $Ug \subset \mathcal{L}(\xi_i)$ and hence $\mathcal{L}(\xi_i)$ is open in $\mathcal{L}$.

Let $\{g_n\}_{n = 1}^\infty$ be a sequence in $\mathcal{L}(\xi_i)$ and suppose $g_n \rightarrow g \in \mathcal{L}$ for $n \rightarrow \infty$. Then
\begin{align*}
d^\perp(\xi_i(g \cdot p), dg (\xi_i(p)))
& \leq
d^\perp(\xi_i(g \cdot p), \xi_i(g_n \cdot p)) + d^\perp (dg_n (\xi_i(p)), dg (\xi_i(p))),
\end{align*}
where $d^\perp$ denotes the Riemannian distance on $T^\perp N$. The topology induced by $d^\perp$ coincides with the original topology on $T^\perp N$ \cite[Proposition 6.1 in Chapter VII]{Lan99}. Thus by continuity of $\xi_i : N \rightarrow T^\perp N$ and of $\mathcal{L} \rightarrow T^\perp N$, $g \mapsto dg(\xi_i(p))$, the right term of the above inequality converges to zero. Hence $\xi_i(g \cdot p) = dg (\xi_i(p))$ and therefore $\mathcal{L}(\xi_i)$ is closed in $\mathcal{L}$.

Since $\mathcal{L}$ is connected and $\mathcal{L}(\xi_i)$ is not empty, we have $\mathcal{L}(\xi_i) = \mathcal{L}$. This implies that the slice representation of $N$ is trivial. Thus $N$ is a principal orbit.
\end{proof}

\begin{remark}\label{rembco}
If the group $\mathcal{L}$ is not connected, there exists a counterexample to the theorem. A finite dimensional counterexample is given by the action of $\mathbb{R} \times \{\pm 1\}$ on $\mathbb{E}^2$ defined by $(a, \pm1) \cdot (x,y) := (x + a, \pm y) $; the orbit through the origin is exceptional. For an infinite dimensional counterexample, see Remark \ref{rem5}. 
\end{remark}


\begin{remark} \label{rembco2}
Similarly to the proof of Theorem \ref{thm2}, one can solve \cite[Exercise 2.11.19]{BCO}. Note that here the group $G$ should be connected.
\end{remark}

\section{Proof of Corollary \ref{thm1}} \label{secI}

Let $G$ be a connected compact semisimple Lie group and $K$ a symmetric subgroup of $G$, that is, there exists an involutive automorphism $\theta$ of $G$ satisfying $G_\theta^0 \subset K \subset G_\theta$, where $G_\theta$ denotes the fixed point subgroup of $\theta$ and $G_\theta^0$ its identity component. The differential of $\theta$ is still denoted by $\theta$. Denote by $\mathfrak{g} = \mathfrak{k} + \mathfrak{m}$ the decomposition into the $(\pm1)$-eigenspaces of $\theta$, which is called the \emph{canonical decomposition} associated to $(\mathfrak{g}, \mathfrak{k})$. Assume that $\mathfrak{g}$ is equipped with an inner product which is invariant by $\operatorname{Ad}(G)$ and $\theta$. Then the canonical decomposition $\mathfrak{g} = \mathfrak{k} + \mathfrak{m}$ is orthogonal. We equip $G$ with the corresponding bi-invariant Riemannian metric and $G/K$ with the normal homogeneous metric. Then $G/K$ is a symmetric space of compact type and the projection $\pi: G \rightarrow G/K$ is a Riemannian submersion.

Let $H$ be a closed subgroup of $G$, which acts on $G/K$ isometrically by $b \cdot aK := (ba) K$. We know that the $H$-action on $G/K$ is hyperpolar if and only if the $(H \times K)$-action on $G$ defined by $(b, c) \cdot a = bac^{-1}$ is hyperpolar \cite[Proposition 2.11]{HPTT95}. It is easy to see that the orbit $(H \times K) \cdot a$ is the inverse image of the orbit $H \cdot aK$ under $\pi$. In particular those orbits have the same codimension. It is also easy to see that $(H \times K) \cdot a$ is principal if and only if $H \cdot aK$ is principal. Thus $(H \times K) \cdot a$ is exceptional if and only if $H \cdot aK$ is exceptional. 

We denote by $\mathcal{G} := H^1([0,1], G)$ the Hilbert Lie group of all Sobolev $H^1$-paths from $[0, 1]$ to $G$ and by $V_\mathfrak{g} := H^0([0,1], \mathfrak{g})$ the Hilbert space of all $H^0$-maps from $[0,1]$ to the Lie algebra $\mathfrak{g}$ of $G$. Consider the subgroup
\begin{equation*}
P(G, L) := \{g \in \mathcal{G} \mid (g(0), g(1)) \in L\}
\end{equation*}
for a closed subgroup $L$ of $G \times G$. Then $P(G, L)$ acts on $V_\mathfrak{g}$ by the affine isometry:
\begin{equation*}
g * u = gug^{-1} - g' g^{-1}.
\end{equation*}
This action is isometric and PF \cite[p.\ 132]{Ter95}. 

The \emph{parallel transport map} $\Phi: V_\mathfrak{g} \rightarrow G$ is defined by $\Phi(u) := g_u(1)$ where $g_u \in \mathcal{G}$ denotes the unique solution to the linear ordinary differential equation $g^{-1} g' = u$, $g(0) = e$. It is easy to see that the restriction of $\Phi$ to the space of constant paths $\hat{\mathfrak{g}}$ with values in $\mathfrak{g}$ is identified with the exponential map $\exp :\mathfrak{g} \rightarrow G$. In particular $\Phi$ is surjective. Moreover (\cite[Proposition 1.1]{Ter95}, \cite[Theorem 4.5]{TT95}):


\begin{proposition}[Terng \cite{Ter95}, Terng-Thorbergsson \cite{TT95}]\ 
\begin{enumerate}
\setlength{\itemsep}{1.5mm}
\item $\Phi$ is a Riemannian submersion,
\item $\Phi(g * u) = g(0) \Phi(u) g(1)^{-1}$ for $g \in \mathcal{G}$ and $u \in V_\mathfrak{g}$,
\item $P(G, L) * u = \Phi^{-1}(L \cdot \Phi(u))$, where $L$ acts on $G$ by $(b, c) \cdot a := bac^{-1}$.
\end{enumerate}
\end{proposition}


In particular, the orbits $P(G, L) * u$ and $L \cdot \Phi(u)$ have the same codimension. Moreover \cite[Theorem 1.2]{Ter95}: 


\begin{proposition}[Terng \cite{Ter95}]\ \label{propterng}
Suppose that the $L$-action on $G$ is hyperpolar. Then
\begin{enumerate}
\setlength{\itemsep}{1.5mm}
\item the $P(G, L)$-action on $V_\mathfrak{g}$ is polar, 
\item the orbit $P(G, L) * u$ is principal if and only if the orbit $L \cdot \Phi(u)$ is principal.
\end{enumerate}
\end{proposition}


In particular, $P(G, L) * u$ is exceptional if and only if $L \cdot \Phi(u)$ is exceptional. In order to apply Theorem \ref{thm2} to the $P(G, L)$-action, we need:


\begin{lemma}\label{lemconnected}
Suppose that $G$ is simply connected and $L$ is connected. Then $P(G, L)$ is connected.
\end{lemma}
\begin{proof}
Consider the map $\Psi: \mathcal{G} \rightarrow G \times G$ defined by $\Psi(g) = (g(0), g(1))$, which is a submersion \cite[p.\ 132]{Ter95}. It is easy to see that the based loop group $P(G, \{e\} \times \{e\})$ acts on each fiber of $\mathcal{G}$ simply transitively. Consider the restriction
\begin{equation*}
\Psi: P(G, L) \rightarrow L.
\end{equation*}
We have the exact sequence of homotopy groups
\begin{equation*}
\cdots
\rightarrow
\pi_0(P(G, \{e\} \times \{e\}))
\rightarrow
\pi_0(P(G, L))
\rightarrow
\pi_0(L).
\end{equation*}
Since $G$ is simply connected, $P(G, \{e\} \times \{e\})$ is connected. Thus the sequence shows that $P(G, L)$ is  connected. 
\end{proof}

\begin{proof}[\textbf{\textup{Proof of Corollary \ref{thm1}}}]
From Proposition \ref{propterng}, the $P(G, H \times K)$-action is polar. From Lemma \ref{lemconnected}, $P(G, H \times K)$ is connected. Thus by Theorem \ref{thm2}, the $P(G, H \times K)$-action does not have exceptional orbits. This condition is equivalent to the condition that the $H$-action on $G/K$ does not have exceptional orbits, as discussed above. This proves Corollary \ref{thm1}. 
\end{proof}

\begin{remark}\label{rem5}
If $G/K$ is not simply connected, there exists a counterexample to Corollary \ref{thm1}. For example, the action of $SO(2)$ on the real projective plane $\mathbb{R} P^2 = SO(3)/O(2)$ by rotation has an exceptional orbit (equator). In this case, from the above discussions, the corresponding $P(G, H \times K)$-orbit is also exceptional and the group $P(G, H \times K)$ is not connected since $K$ is not connected. This gives an infinite dimensional counterexample to Theorem \ref{thm2} when $\mathcal{L}$ is not connected.
\end{remark}

\section*{Acknowledgements}
The author would like to thank Professor Yoshihiro Ohnita for many useful discussions and invaluable suggestions. The author is also grateful to Professor Takashi Sakai for his constant support.


\end{document}